\RequirePackage{fix-cm}
\documentclass[smallextended]{svjour3}       
\smartqed  
\usepackage{graphicx}
\usepackage{amsmath}
\numberwithin{equation}{section}
\numberwithin{corollary}{section}

\numberwithin{definition}{section}
\numberwithin{proposition}{section}
\numberwithin{lemma}{section}
\begin{document}

\title{Generalization of Lax Equivalence Theorem on Unbounded Self-adjoint Operators with Applications to Schr\"{o}dinger Operators 
}


\author{Yidong Luo* 
}


\institute{Yidong Luo \at
              School of Mathematics and Statistics, Wuhan University, Province Hubei, China 430072 \\
              \email{Sylois@whu.edu.cn}           
}

\date{Received: date / Accepted: date}

\maketitle

\begin{abstract}
Define $ A $ a unbounded self-adjoint operator on Hilbert space $ X $. Let $ \{ A_n \} $ be its resolvent approximation sequence with closed range $ \mathcal{R}(A_n) (n \in \mathrm{N}) $, that is, $ A_n (n \in \mathrm{N}) $ are all self-adjoint on Hilbert space $ X $ and
\begin{equation*}
\hbox{ \raise-2mm\hbox{$\textstyle s-\lim \atop \scriptstyle {n \to \infty}$}} R_\lambda (A_n)
 = R_\lambda (A)\quad  (\lambda \in \mathrm{C} \setminus \mathrm{R}), \ \textrm{where}  \ R_ \lambda(A) := (\lambda I-A)^{-1}.
\end{equation*}
The Moore-Penrose inverse $ A^\dagger_n \in \mathcal{B}(X) $ is a natural approximation to the Moore-Penrose inverse $ A^\dagger $.
This paper shows that: $ A^\dagger $ is continuous and strongly converged by $ \{ A^\dagger_n \} $ if and only if
$ \sup\limits_n \Vert A^\dagger_n \Vert < +\infty $.
On the other hand, this result tells that arbitrary bounded computational scheme $ \{ A^\dagger_n \} $ induced by resolvent approximation $ \{ A_n \} $ is naturally instable (that is, $ \sup_n \Vert A^\dagger_n \Vert = \infty $) for any self-adjoint operator equation with non-closed range, for example, free Schr\"{o}dinger operator, Schr\"{o}dinger operator with Coulumb potential and Schr\"{o}dinger operator in model of many particles. This implies the infeasibility to globally and approximately solve non-closed range self-ajoint operator equation by resolvent approximation.
\end{abstract}

\section{Introduction}
In practical fields, a class of problems can be formulated as operator equation
\begin{equation}
Ax = y
\end{equation}
where $ A $, in general, is a bounded operator mapping from Hilbert space $ X $ to $ Y $.
In many cases, directly solving (1.1) is infeasible or the solution of (1.1) is severely sensitive on the R.H.S $ y $. It is common in numerical analysis to solve a system
\begin{equation}
A_n x = y
\end{equation}
where $ A_n: X \to Y $ is a bounded operator approximate to $ A $ in some sense such that (1.2) is easier to solve and the solution of (1.2) depends stably on the R.H.S $ y $. Naturally, one will wonder how to
design the operator $ A_n $ in order that solution of (1.2) can be a good approximation to solution of (1.1).
\newline \indent In the ideal case when $ A, A_n $ are all bijective, the classical Lax equivalence theorem answer to this question, and illustrates that
\begin{theorem}
If $ A_n $ is approximate to $ A $ in the following sense
\begin{equation*}
(\text{consistency}): \ \ \Vert A_n x - A x \Vert \to 0 \quad  (x \in X)
\end{equation*}
 then it follows that
 \begin{equation*}
 (\text{convergence}): \ \ \Vert A^{-1}_n x - A^{-1} x \Vert \to 0  \quad (x \in Y), \Longleftrightarrow (\text{stability}): \ \sup_n \Vert A^{-1}_n \Vert <\infty.
 \end{equation*}
\end{theorem}
\indent This result tells us theoretically that, if a consistent approximation sequence $ \{ A_n \} $ possesses stability, then it is properly designed. Also notice that when $ A, A_n $ are all bijective, then the inverses of them will also be bounded from $ Y $ to $ X $. Thus the sensitivity problem on the R.H.S $ y $ vanishes simultaneously.
\newline \indent Further, it is a natural idea to establish similar criteria for a wider class of linear operators. However, in general case one will always face a linear operator (not necessarily bounded) which is neither injective ($ \mathcal{N}(A) \neq 0 $) nor surjective ($ \mathcal{R}(A) \neq Y $). Now $ A^{-1} $ does not exist, it is necessary to introduce the generalized inverse $ A^\dagger $ to develop the extension result of Theorem 1.1.
 \newline \textbf{Moore-Penrose inverse of Linear Operators}: Let $ X,Y $ be Hilbert spaces. For linear operator $ A: \mathcal{D}(A) \subseteq X \longrightarrow Y $,  we denote $ \mathcal{D}(A) $, $ \mathcal{R}(A) $, $\mathcal{N}(A) $, $\mathcal{G}(A) $ as its domain, range, kernel and graph respectively. If its domain is decomposable with respect to the kernel space, that is,
\begin{equation}
\mathcal{D}(A) = \mathcal{N}(A) \oplus \mathcal{C}(A), \ \textrm{where} \ \mathcal{C}(A) := \mathcal{D}(A) \cap \mathcal{N}(A)^\perp,
\end{equation}
then we can define $ A_0 := A|_{\mathcal{C}(A)} $ and $ A^{-1}_0: \mathcal{R}(A) \subseteq Y \longrightarrow \mathcal{C}(A) \subseteq X $ exists. Regard $ A^{-1}_0 $ as $ A^\dagger |_{\mathcal{R}(A)} $ and extend it  to $ A^\dagger $ with
 \begin{equation}
 \mathcal{D}(A^\dagger) = \mathcal{R}(A)+\mathcal{R}(A)^\perp,
\end{equation}
\begin{equation}
 \mathcal{N}(A^\dagger) = \mathcal{R}(A)^\perp.
 \end{equation}
 Above extension is unique and well defined. This defines the Moore-Penrose inverse $ A^\dagger $ of linear operator $ A $ (also denoted as the maximal Tseng inverse, see [1,Chapter 9.3, Definition 2]).
 \newline \indent Notice that, if $ A $ is closed, then $\mathcal{N}(A) $ is closed, and recall the fact that a space in Hilbert space is decomposable with respect to any closed subspace (See [1, Chapter 9.2, Ex.5]), then (1.3) automatically holds. Thus, $ A^\dagger $ exists throughout this paper. Moreover, by [1, Chapter 9.3, Ex13],
 \begin{equation}
  A^\dagger \in \mathcal{B}(Y,X) \Longleftrightarrow  \mathcal{R}(A) \ \textrm{closed} \Longleftrightarrow \mathcal{D}(A^\dagger) = Y.
 \end{equation}
  For more about the Moore-Penrose inverse, see [1, Chapter 9].
 \begin{remark}
 $ A^\dagger $ helps give a more generalized definition for "solution": for $ A \in \mathcal{C}(X,Y) $ neither injective nor surjective with not necessarily closed range $ \mathcal{R}(A) $:
 \newline Case I: when $ y \in \mathcal{R}(A) $, (1.1) has infinitely many solutions. $ A^\dagger y $ gives the minimun norm solution for (1.1) which can be convenient for the preceding discussion on convergence.
 \newline Case II: when $ y \in \mathcal{R}(A) + \mathcal{R}(A)^\perp, P_{\mathcal{R}(A)^\perp} y \neq 0 $, (1.1) has no solution. But $ u:= A^\dagger y $  provides the best approximate solution in the sense that
 \begin{equation*}
 \Vert A u - y \Vert = \inf_{x \in \mathcal{C}(A)} \Vert Ax - y \Vert.
 \end{equation*}
 Case III: when $ y \in Y \setminus (\mathcal{R}(A) + \mathcal{R}(A)^\perp) $, (1.1) has no solution, and this part is also not attainable for $ A^\dagger $. But when the range space is closed, this case will simultaneously vanish.
 \newline \indent Thus we see that we can construct a more complete Consistency-Stability-Convergence framework facing a wider class of R.H.S. $ y $ with Moore-Penrose inverse.
 \end{remark}
\indent
Let $  \mathcal{L}(X,Y) $ denote the set of  all linear operators mapping from $  X $ to $ Y $, $ \mathcal{C}(X,Y) $ the set of all $ A \in \mathcal{L}(X,Y) $ with closed graph, $ \mathcal{B}(X,Y) $ the set of  all bounded linear operators $ A \in \mathcal{L}(X,Y)$, and $ \mathcal{CR}(X,Y) $ the set of  all $ A \in \mathcal{B}(X,Y)$ with closed range. When it concerns operator $ A \in \mathcal{B}(X,Y) $, we assume that $ \mathcal{D}(A) = X $.
\newline \indent We recall the definitions of consistency, stability and convergence (refer to [10]):
\newline (A1): Strong consistency:
\begin{equation*}
 \Vert A_n x - A x \Vert \to 0 \ ( n \to \infty) \quad (x \in X).
\end{equation*}
(A2): Uniform consistency:
\begin{equation*}
\Vert A_n - A \Vert \to 0 \ ( n \to \infty).
\end{equation*}
 (B1): Stability:
\begin{equation*}
\sup_n \Vert A^\dagger_n \Vert <\infty.
\end{equation*}
(C1): Perfect strong convergence:
\begin{equation*}
\mathcal{D}(A^\dagger)=Y, \hbox{ \raise-2mm\hbox{$\textstyle s-\lim \atop \scriptstyle {n \to \infty}$}} A^\dagger_n=A^\dagger.
\end{equation*}
(C2): Perfect uniform convergence:
\begin{equation*}
 \mathcal{D}(A^\dagger)=Y, \hbox{ \raise-2mm\hbox{$\textstyle \lim \atop \scriptstyle {n \to \infty}$}} \Vert A^\dagger_n - A^\dagger \Vert = 0.
\end{equation*}
\begin{remark}
The formulations (A1)-(C2) are originally defined in [10,11] for  $ A_n, A \in \mathcal{B}(X,Y) \ (n \in \mathrm{N})$. In this paper we will propose a new type of consistency to replace (A1) and (A2) for unbounded operator $ A \in \mathcal{C}(X,Y) $ and its approximation operator sequence $ \{ A_n \} \subseteq  \mathcal{C}(X,Y) $ with closed range $ \mathcal{R}(A_n) $. As to (B1)-(C2), they are retained and proven well-defined for above setting: By (1.4), it yields that $ A^\dagger_n \in \mathcal{B}(Y,X) (n \in \mathrm{N})$. Thus $ \Vert A^\dagger_n \Vert (n \in \mathrm{N}) $ are all finite and $ (B1) $ is well defined. Besides, provided with the original operator $ A \in \mathcal{C}(X,Y) $ and $ \mathcal{D}(A^\dagger)=Y $, by (1.4) we have $ A^\dagger \in \mathcal{B}(Y,X) $. Thus we can discuss the strong convergence and the norm convergence of $ A^\dagger_n $ to $ A^\dagger $ in sense of (C1) and (C2).
\end{remark}
 \textbf{\textrm{Previous results and main result}}: For $ A \in \mathcal{CR}(X,Y) $, provided with
approximation sequence $ \{ A_n \} $ in  $ \mathcal{CR}(X,Y) $,
it is shown in [14] that,
\begin{itemize}
\item
If $ \{ A_n \} $ and $ A $ satisfy $ (A2) $, then $ (C2) \Longleftrightarrow  (B1) $;
\item
If $ \{ A_n \} $ and $ A $ satisfy $ (A1) $, then
\begin{equation}
(C1) \Longleftrightarrow
(B1) \ \textrm{and} \
\begin{cases}
A^\dagger_n A_n \stackrel{s}\longrightarrow A^\dagger A \\
A_n  A^\dagger_n \stackrel{s}\longrightarrow A A^\dagger.
\end{cases}
 \end{equation}
 \end{itemize}

\indent Above results are all based on a priori information that $ A $ possesses a closed range.
Without this assumption, some improved versions of above results are given for  $ A \in \mathcal{B}(X,Y) $ in [10].
 \newline \indent For $ A \in \mathcal{B}(X,Y) $, provided with
approximation sequence $ \{ A_n \} $ in  $ \mathcal{CR}(X,Y) $,
\begin{itemize}
\item
If $ \{ A_n \} $ and $ A $ satisfy $ (A2) $,
 then (B1) implies that $ A $ possesses the closed range $ \mathcal{R}(A) $. Furthermore, $ (B1) \Longleftrightarrow  (C1) \Longleftrightarrow (C2) $.
 \item
 If $ \{ A_n \} $ and $ A $ satisfy $ (A1) $, then
\begin{equation}
(C1)
 \Longleftrightarrow
(B1) \ \textrm{and} \
\begin{cases}
\hbox{ \raise-2mm\hbox{$\textstyle s-\lim \atop \scriptstyle {n \to \infty}$}} \ \mathcal{R}(A_n)
=
\hbox{ \raise-2mm\hbox{$\textstyle  w-\widetilde{\lim} \atop \scriptstyle {n \to \infty}$}} \mathcal{R}(A_n)
=\mathcal{R}(A) \\
\hbox{ \raise-2mm\hbox{$\textstyle s-\lim \atop \scriptstyle {n \to \infty}$}} \ \mathcal{N}(A_n)
=
\hbox{ \raise-2mm\hbox{$\textstyle  w-\widetilde{\lim} \atop \scriptstyle {n \to \infty}$}} \mathcal{N}(A_n)
=\mathcal{N}(A).
\end{cases}
 \end{equation}
 \end{itemize}

\indent The equivalences in (1.5) and (1.6) are expressed by (B1) and additional conditions. Eliminating these additional conditions but supplementing self-adjoint assumptions for $ A $ and $ \{ A_n \} $, the equivalence result between (B1) and (C1) (under (A1)) is obtained in [12].
\newline \indent
This paper intends to generalize above result into a unbounded case. Before we formulate the main result, we indicate a type of new consistency since the consistency (A1) and (A2) does not suit the approximation of unbounded operators any more.
 \newline (A3): Resolvent consistency:
 \newline Suppose that $ A $ and $ \{ A_n \} $ are all self-adjoint operators (possible unbounded)
on Hilbert space $ X $. If
\begin{equation*}
\hbox{ \raise-2mm\hbox{$\textstyle s-\lim \atop \scriptstyle {n \to \infty}$}} R_\lambda (A_n)
 = R_\lambda (A)\quad  (\lambda \in \mathrm{C} \setminus \mathrm{R}), \  \ R_ \lambda(A) := (\lambda I-A)^{-1},
\end{equation*}
then we say that $ \{ A_n \} $ and $ A $ satisfy the resolvent consistency, i.e., $ \hbox{ \raise-2mm\hbox{$\textstyle s.r.s-\lim \atop \scriptstyle {n \to \infty}$}} A_n =A $.
\newline \indent Now we give the main result:
\begin{theorem}
 Let \begin{math} A \end{math} be self-adjoint operator(possible unbounded) on Hilbert space \begin{math} X \end{math},
 \begin{math} \{ A_n \} \end{math} a sequence of self-adjoint operators on \begin{math} X \end{math} with closed range \begin{math} \mathcal{R}(A_n) (n \in \mathrm{N}) \end{math}. If \begin{math} \{ A_n \} \end{math} and $ A $ satisfy the resolvent consistency, then
 \newline \indent (a) $ \sup_n \Vert A^\dagger_n \Vert <+\infty \ ((B1)) \Longrightarrow A \text{  preserves closed range } \mathcal{R}(A) $;
  \newline \indent (b) $ \mathcal{D}(A^\dagger)=X,
 A^\dagger_n  \stackrel{s}\longrightarrow  A^\dagger \ ((C1)) \Longleftrightarrow \sup_n \Vert A^\dagger_n \Vert <+\infty \ ((B1)) $.
 \end{theorem}
\indent The remainder of this paper is organized as follows:
 In section 2, we introduce some basic conceptions, such as unbounded self-adjoint operators and the strong graph limit.
 In section 3 and section 4, we prove the results (a) and (b) respectively. In section 5, we give the applications of the main result to Schrodinger operator equation. In section 6, we give some corollaries. In section 7, we conclude the main work of this paper and give prospects of future work.
\section{Preliminary and Basic Lemmas}
\subsection{Moore-Penrose inverse}
\begin{proposition}
For a densely defined closed operator $ A $ on Hilbert space $ X $, its Moore-Penrose inverse $ A^\dagger $
satisfies the following two identities
\begin{equation}
A^\dagger A x = P_{\overline{\mathcal{C}(A)}} x\quad (x \in \mathcal{N} ( A ) \oplus \mathcal{C} (A)),
\end{equation}
\begin{equation}
 A A^\dagger y =P_{\overline{\mathcal{R}(A)}} y \quad (y \in \mathcal{R} ( A ) \oplus \mathcal{R} (A)^{\perp}).
\end{equation}
\end{proposition}
\begin{proof}
This result can be found in [1, Chapter 9.3, Theorem 1]. However, for the convenience of readers, we provide a proof of $ (2.1) $ here, $ (2.2) $ could be obtained in a similar way.
\newline\indent  For $ x \in \mathcal{D}(A) = \mathcal{N}(A) \oplus \mathcal{C}(A) $, it can be uniquely represented as
\begin{equation*}
x = x_1 +x_2, \ \textrm{where} \ x_1 \in \mathcal{N}(A), \ x_2 \in \mathcal{C}(A), \ \textrm{and} \ x_1 \perp x_2.
\end{equation*}
The L. H. S. of $ (2.1) $ reads as follows
\begin{equation*}
A^\dagger A x = A^\dagger A ( x_1 + x_2 ) = A^\dagger A x_2 = A^\dagger A_0 x_2 = A^{-1}_0 A_0 x_2 = x_2 = P_{\overline{\mathcal{C}(A)}} x.
\end{equation*}
\qed
\end{proof}
\subsection{Unbounded self-adjoint operator and the strong graph convergence}
We firstly introduce the concept of adjoint operator.
 \begin{definition}
 Let \begin{math} A \end{math} be a densely defined closed operator on Hilbert space $ X $. Set
\begin{equation*}
 \mathcal{D} (A^*) := \{ u \in X | \  \textrm{There exists a} \ v \in X, \ \text{such that} \
 \langle u, Ax \rangle = \langle v , x \rangle \quad (x \in \mathcal{D}(A)) \}.
 \end{equation*}
Then
\begin{eqnarray*}
A^* :  \mathcal{D}(A^*) \subseteq X & \longrightarrow & X \\
 u & \longmapsto & v
\end{eqnarray*}
is defined as the adjoint operator of $ A $,
 where \begin{math} \mathcal{D}(A^*) \end{math} is the domain of \begin{math} A^* \end{math}.
\end{definition}
\begin{definition}
Let \begin{math} A \end{math} be a densely defined closed linear operator on Hilbert space $ X $.
 If $ A = A^* $, then we call \begin{math} A \end{math} self-adjoint. Notice that $ A = A^* $ means:
\newline \indent \ (1) \ $ \mathcal{D}(A) = \mathcal{D} (A^*) $,
\newline \indent \ (2) \ $  \langle Ax,y \rangle  = \langle x,Ay \rangle \quad (x,y \in \mathcal{D}(A)). $
\end{definition}
For unbounded self-adjoint operator (actually not restricted in this case), we additionally introduce a convergence of new type:
\begin{definition}
 Let $ \{ A_n \} $ be a sequence of  closed linear operators on Hilbert space $ X $. We define
 \begin{equation*}
 \hbox{ \raise-2mm\hbox{$\textstyle s-\lim \atop \scriptstyle {n \to \infty}$}}  \mathcal{G}(A_n) :=
 \end{equation*}
 \begin{equation*}
 \{(u,v) \in X \times X :
 \textrm{There exists a} \ u_n \in \mathcal{D} (A_n) (n \in \mathrm{N}) \ \textrm{such that} \ (u_n,A_n u_n) \stackrel{s}{\longrightarrow} (u,v)  \}
 \end{equation*}
 If $ \hbox{ \raise-2mm\hbox{$\textstyle s-\lim \atop \scriptstyle {n \to \infty}$}}  \mathcal{G}(A_n) $ is the graph of an operator $ A $, then we say that $ A $ is the strong graph limit of $ \{  A_n \} $ and write $ \hbox{ \raise-2mm\hbox{$\textstyle s.g. -\lim \atop \scriptstyle {n \to \infty}$}} A_n=A $.
\end{definition}
\indent The following result indicates that the resolvent convergence and the strong graph convergence are equivalent when $ A_n (n \in \mathrm{N}) $ and $ A $ are all self-adjoint.
\begin{lemma}
Let \begin{math} A_n(n \in \mathrm{N}), \ A \end{math} be self-adjoint operators on Hilbert space \begin{math} X \end{math}, then
\begin{equation*}
\hbox{ \raise-2mm\hbox{$\textstyle s.r.s-\lim \atop \scriptstyle {n \to \infty}$}} A_n =A \Longleftrightarrow \hbox{ \raise-2mm\hbox{$\textstyle s.g-\lim \atop \scriptstyle {n \to \infty}$}} A_n=A.
\end{equation*}
\end{lemma}
\begin{proof}
See [8, P.293 Theorem VIII. 26]. \qed
\end{proof}
\subsection{Characterization for convergence of orthogonal projection sequence}
Let $ \{ X_n \} $ be a subspace sequence of Hilbert space $ X $. We define
 \begin{equation*}
 \hbox{ \raise-2mm\hbox{$\textstyle s-\lim \atop \scriptstyle {n \to \infty}$}} X_n :=\{ x \in X: \textrm{There exists a} \ x_n \in X_n(n \in \mathrm{N}) \ \textrm{such that} \  x_n \stackrel{s}{\longrightarrow} \  x \}
 \end{equation*}
 and
 \begin{equation*}
 \hbox{ \raise-2mm\hbox{$\textstyle  w-\widetilde{\lim} \atop \scriptstyle {n \to \infty}$}} X_n :=\{ x \in X: \textrm{There exists a} \ x_n \in X_{k_n}(n \in \mathrm{N}) \ \textrm{such that} \ x_n\stackrel{w}{\longrightarrow}  x \}.
 \end{equation*}
\indent The convergence of orthogonal projection sequence $ \{ P_{X_n} \} $ is characterized in the following result.
\begin{lemma}
Let $ X $ be Hilbert space and $ \{ X_n \} $ a sequence of closed subspaces of $ X $,
Then
\begin{equation*}
\{ P_{X_n} \} \ \text{is strongly convergent }
\Longleftrightarrow
\hbox{ \raise-2mm\hbox{$\textstyle s-\lim \atop \scriptstyle {n \to \infty}$}} X_n
=
\hbox{ \raise-2mm\hbox{$\textstyle  w-\widetilde{\lim} \atop \scriptstyle {n \to \infty}$}} X_n.
\end{equation*}
Moreover, in the case that \begin{math} \{ P_{X_n} \} \end{math}  is strongly convergent,
\begin{equation*}
 \hbox{ \raise-2mm\hbox{$\textstyle s-\lim \atop \scriptstyle {n \to \infty}$}} P_{X_n} =P_M , \text{where} \  M:= \hbox{ \raise-2mm\hbox{$\textstyle s-\lim \atop \scriptstyle {n \to \infty}$}} X_n.
 \end{equation*}
\end{lemma}
\begin{proof}
See [10, Lemma 2.13]. \qed
\end{proof}
\subsection{Weak convergence}
\begin{lemma}
Let \begin{math} X \end{math} be a Hilbert space, \begin{math} \{ x_n \} \end{math} a weakly convergent sequence of \begin{math} X \end{math} with \begin{math} x_{\infty}= \hbox{ \raise-2mm\hbox{$\textstyle w-\lim \atop \scriptstyle {n \to \infty}$}} x_n \end{math}. Then
\begin{equation*}
\sup\limits_n \Vert x_n \Vert < +\infty,\
\Vert x_{\infty} \Vert \leq \hbox{ \raise-2mm\hbox{$\textstyle \underline{\lim} \atop \scriptstyle {n \to \infty}$}} \Vert x_n \Vert.
\end{equation*}
\end{lemma}
\begin{proof}
See [7, p.120, Theorem 1]. \qed
\end{proof}
\section{Proof of Result (a)}
Before the proof for result (a), we first prepare two lemmas to describe how the kernel space sequence $ \{ \mathcal{N}(A_n) \} $
converges in a strong and weak sense.
\begin{lemma}
 Let $ A $ be a closed linear operator, $ \{ A_n \} $ a sequence
 of closed linear operators with closed range $ \mathcal{R}(A_n)(n \in \mathrm{N}) $.
 Suppose
\begin{equation}
\hbox{ \raise-2mm\hbox{$\textstyle s.g-\lim \atop \scriptstyle {n \to \infty}$}} A_n =A
\end{equation}
and
\begin{displaymath}
 \sup\limits_n\Vert A^\dagger_n \Vert <\infty.
 \end{displaymath}
  Then, for $ y \in \mathcal{R}(A) $ and any sequence $ \{ y_n \}$ such that
  \begin{displaymath}
  y_n \in \mathcal{R}(A_n) \ \textrm{and} \  y_n \stackrel{s}{\longrightarrow} y.
  \end{displaymath}
  We have
\begin{displaymath}
\hbox{ \raise-2mm\hbox{$\textstyle s-\lim \atop \scriptstyle {n \to \infty}$}} A_n^{-1} (y_n) =A^{-1} (y) .
\end{displaymath}
Denote that $ A^{-1}(y)=\{x \in \mathcal{D}(A): Ax=y \} $. Furthermore, setting $ y_n = y = 0 \ (n \in \mathrm{N}) $, it follows that
\begin{displaymath}
 \hbox{ \raise-2mm\hbox{$\textstyle s-\lim \atop \scriptstyle {n \to \infty}$}} \mathcal{N}(A_n) = \mathcal{N}(A).
 \end{displaymath}
\end{lemma}
\begin{proof}
  Let $ y \in \mathcal{R}(A) $ and $ \{ y_n \}$ be any sequence such that
  \begin{displaymath}
  y_n \in \mathcal{R}(A_n) (n \in \mathrm{N}) \ \textrm{and} \  y_n \stackrel{s}{\longrightarrow} y.
  \end{displaymath}
  \indent "$ \subseteq $": Suppose that \begin{math} x \in \hbox{ \raise-2mm\hbox{$\textstyle s-\lim \atop \scriptstyle {n \to \infty}$}} A^{-1}_n(y_n) \end{math}. There exist a sequence \begin{math} \{ x_n \} \end{math} such that
 \begin{displaymath}
 x_n \in A^{-1}_n(y_n) ( \forall n \in \mathrm{N}) \  \textrm{and} \  x_n  \stackrel{s}{\longrightarrow} x.
 \end{displaymath}
Notice that
 \begin{equation*}
 A_n x_n = y_n  \stackrel{s}{\longrightarrow} y.
  \end{equation*}
 We have
  \begin{displaymath}
  (x_n, A_n x_n) \stackrel{s}{\longrightarrow} (x,y) \textrm{ in } X \times X.
  \end{displaymath}
  Since $ \hbox{ \raise-2mm\hbox{$\textstyle s-\lim \atop \scriptstyle {n \to \infty}$}} \mathcal{G}(A_n) =\mathcal{G}(A) \ (\textrm{by} \ (3.1) )$,  we have $(x,y) \in \mathcal{G}(A)$, that is, $ x \in \mathcal{D}(A), y=Ax $. So \begin{math} x \in A^{-1}(y) \end{math}.
  \newline \indent "$\supseteq$": Assume that $x \in A^{-1}(y)$. Then
  \begin{displaymath}
 (x,y) \in \mathcal{G}(A) \stackrel{(3.1)}{=} \hbox{ \raise-2mm\hbox{$\textstyle s-\lim \atop \scriptstyle {n \to \infty}$}} \mathcal{G}(A_n).
  \end{displaymath}
There exists a sequence $ (x_n, A_n x_n) \in \mathcal{G}(A_n) $ such that
  \begin{equation}
(x_n, A_n x_n) \stackrel{s}{\longrightarrow} (x,y).
  \end{equation}
  In the following, we set
  \begin{displaymath}
   u_n =A^\dagger_n(y_n-A_n(x_n)), \ p_n := x_n + u_n \ (n \in \mathrm{N})
   \end{displaymath}
  and prove
  \begin{displaymath}
  p_n \in A^{-1}_n(y_n), \ p_n \stackrel{s}{\longrightarrow} x.
  \end{displaymath}
\indent First, we can check that
 \begin{displaymath}
 \Vert p_n -x \Vert = \Vert u_n+x_n-x \Vert
 \leq \Vert A^\dagger _n (y_n-A_n x_n)\Vert +\Vert x_n -x \Vert
  \end{displaymath}
\begin{displaymath}
\leq M \Vert y_n -A_n x_n\Vert +\Vert x_n -x \Vert \stackrel{n \to \infty}{\longrightarrow} 0 \ (\textrm{by} \ (3.2) \ \textrm{and} \ (B1))
\end{displaymath}
where $ M:=\sup_n \Vert A^\dagger_n \Vert.$
\newline \indent Second, for $ A_n u_n = A_n A^\dagger_n(y_n-A_n(x_n)) $, using $ y_n -A_n x_n \in \mathcal{R}(A_n) $ and $ (2.2) $, we have $ A_n u_n = y_n-A_nx_n $. Hence $ A_n p_n = A_n x_n + A_n u_n = y_n $.
\newline \indent Thus, $ x \in \hbox{ \raise-2mm\hbox{$\textstyle s-\lim \atop \scriptstyle {n \to \infty}$}} A^{-1}_n(y_n).$ \qed
\end{proof}
\begin{lemma}
Let $ A, A_n\quad (n \in \mathrm{N}) $ be self-adjoint operators (possible unbounded) on Hilbert space $ X $.
 If $\hbox{ \raise-2mm\hbox{$\textstyle s.g-\lim \atop \scriptstyle {n \to \infty}$}} A_n =A $, then
 \begin{displaymath}
\hbox{ \raise-2mm\hbox{$\textstyle w-\widetilde{\lim}\atop \scriptstyle {n \to \infty}$}}\mathcal{N}(A_n) \subset \mathcal{N}(A).
\end{displaymath}
\end{lemma}
\begin{proof}
Let $ x \in  \hbox{ \raise-2mm\hbox{$\textstyle w-\widetilde{\lim}\atop \scriptstyle {n \to \infty}$}} \mathcal{N}(A_n) $.
There exists a sequence $\{ x_n \} $ such that
\begin{equation}
x_n \in \mathcal{N}(A_{k_n})(k_n \geq n)
\end{equation}
and
\begin{equation}
 x_n \stackrel{w}{\longrightarrow} x  (n \to \infty) .
\end{equation}
For the proof of $ x \in \mathcal{N}(A) $, it is sufficient to prove
\begin{equation}
 \langle x, Au \rangle = 0 \quad (u \in \mathcal{D}(A)).
 \end{equation}
Since $ A $ is self-adjoint, for any $ u \in \mathcal{D}(A) $,
\begin{equation}
\langle x, Au\rangle = \langle x-x_n,Au\rangle + \langle x_n, Au \rangle \quad (n \in \mathrm{N}),
\end{equation}
where
\begin{equation}
 \langle x-x_n,Au\rangle \to 0 \ (n \to \infty) \ (\textrm{by} \ (3.4)).
 \end{equation}
As to the latter term of R.H.S.,
by (3.1), for $ u \in \mathcal{D}(A) $, there exists $ u_{k_n} \in \mathcal{D}(A_{k_n}) $ such that
\begin{equation}
 (u_{k_n}, A_{k_n} u_{k_n})  \stackrel{s}{\longrightarrow}  (u,Au).
\end{equation}
Notice that,
\begin{eqnarray}
 \langle x_n, Au \rangle & = \langle x_n, Au\rangle - \langle A_{k_n} x_n , u_{k_n} \rangle \  (x_n\in \mathcal{N}(A_{k_n}))  \nonumber \\
                         &  = \langle x_n, Au\rangle - \langle x_n , A_{k_n} u_{k_n} \rangle \nonumber \\
                         &   = \langle x_n, Au-A_{k_n} u_{k_n}\rangle.
\end{eqnarray}
By (3.4) and Lemma 2.3,
\begin{displaymath}
\sup\limits_n \Vert x_n \Vert <+\infty.
\end{displaymath}
Thus, with (3.8) and (3.9), it yields that
\begin{equation}
 \vert \langle x_n, Au \rangle \vert
 \leq \Vert x_n \Vert \Vert Au-A_{k_n}u_{k_n} \Vert \longrightarrow 0 \ (n \to \infty).
\end{equation}
\indent Assuming $ n \to \infty $ in (3.6), and using (3.7) and (3.10), we have (3.5). Thus, $ x \in \mathcal{R}(A)^\perp $.
\newline \indent Since $ \mathcal{R}(A)^\perp = \mathcal{N}(A^*) $ holds for all densely defined operator $ A $ on Hilbert space $ X $ (See [5, Chapter X, Proposition 1.13]) and $ A $ is self-adjoint, we have $ x \in \mathcal{N}(A). $
\newline \indent This completes the proof. \qed
\end{proof}

\noindent \textbf{Proof of Result (a)} \indent This proof follows the main idea of [10, Theorem 2.1 (2.22)]. Throughout the whole proof, we will proceed with setting $ \hbox{ \raise-2mm\hbox{$\textstyle s.g-\lim \atop \scriptstyle {n \to \infty}$}} A_n=A $, that is,
\begin{equation}
\mathcal{G}(A)=\hbox{ \raise-2mm\hbox{$\textstyle s-\lim \atop \scriptstyle {n \to \infty}$}} \mathcal{G}(A_n).
\end{equation}
Since $ A, A_n(n \in \mathrm{N}) $ are all self-adjoint satisfying (A3) and by Lemma 2.1, we have
\begin{displaymath}
\hbox{ \raise-2mm\hbox{$\textstyle s.r.s-\lim \atop \scriptstyle {n \to \infty}$}} A_n =A \Longrightarrow \hbox{ \raise-2mm\hbox{$\textstyle s.g-\lim \atop \scriptstyle {n \to \infty}$}} A_n=A.
\end{displaymath}
\indent Let $ \{ y^{(m)} \} \subseteq \mathcal{R}(A)$ and $ \hbox{ \raise-2mm\hbox{$\textstyle s-\lim \atop \scriptstyle {m \to \infty}$}} y^{(m)} =y $.
\newline
 \textbf{Part I} \indent Construct a sequence of pairs $ \{ (x^{(m)}, y^{(m)}) \} \subseteq \mathcal{G}(A) $ with $ \{ x^{(m)} \} $ bounded.
 We proceed with the following three steps.
\newline \textbf{(1)}: The construction of $ \{ x^{(m)} \} $.
\newline \indent We claim that $ (A^\dagger y^{(m)}, y^{(m)}) \in \mathcal{G}(A) $, since
\begin{displaymath}
 A A^\dagger y^{(m)}\stackrel{(2.2)}{=} P_{\overline{\mathcal{R}(A)}} y^{(m)} = y^{(m)} \quad (m \in \mathrm{N}).
 \end{displaymath}
By (3.11), for every $ m \in \mathrm{N} $, there exists a sequence
\begin{displaymath}
   (x^{(m)}_n, y^{(m)}_n) \in \mathcal{G}(A_n), \ n \in \mathrm{N}
   \end{displaymath}
that is, $ x^{(m)}_n \in \mathcal{D}(A_n), y^{(m)}_n=A_n (x^{(m)}_n) \quad (n \in \mathrm{N}) $,
such that
\begin{equation}
x^{(m)}_n \stackrel{s}{\longrightarrow} A^\dagger(y^{(m)}), \ \ y^{(m)}_n \stackrel{s}{\longrightarrow} y^{(m)} \ (n \to \infty).
\end{equation}
Notice that, with (3.12) and (B1),
\begin{equation}
\sup\limits_n \Vert A^\dagger_n(y^{(m)}_n) \Vert \leq \sup\limits_n \Vert A^\dagger_n \Vert \sup\limits_n \Vert y^{(m)}_n \Vert <\infty.
\end{equation}
 Because of (3.13) and the reflexive property of Hilbert space $ X $, by Eberlein-Shmulyan theorem,
$ \{ A^\dagger_n(y^{(m)}_n) \}^{\infty}_{n=1} $ contains a  weakly convergent subsequence
$\{ A^\dagger_{n_j}(y^{(m)}_{n_j}) \}^{\infty}_{j=1} $. Set
\begin{displaymath}
 x^{(m)}:= \hbox{ \raise-2mm\hbox{$\textstyle w-\lim \atop \scriptstyle {j \to \infty}$}}  A^\dagger_{n_j}(y^{(m)}_{n_j}).
\end{displaymath}
\textbf{(2)}: The proof of $ x^{(m)} \in A^{-1}(y^{(m)})$. That is, $ x^{(m)} \in \mathcal{D}(A), Ax^{(m)} = y^{(m)}$.
 \newline \indent For every \begin{math} m \in \mathrm{N} \end{math}, by (3.12),
 \begin{equation}
 A^\dagger(y^{(m)})-x^{(m)} = \hbox{ \raise-2mm\hbox{$\textstyle w-\lim \atop \scriptstyle {j \to \infty}$}} \ x^{(m)}_{n_j}-A^\dagger_{n_j}(y^{(m)}_{n_j}).
\end{equation}
Since
$ x^{(m)}_{n_j} , \ A^\dagger_{n_j} (y^{(m)}_{n_j}) \in A^{-1}_{n_j}(y^{(m)}_{n_j}) $, we can verify \begin{math} x^{(m)}_{n_j}-  A^\dagger_{n_j} (y^{(m)}_{n_j}) \in \mathcal{N}(A_{n_j}) \end{math} for every \begin{math} m \in \mathrm{N} \end{math}. Further by (3.14), we know
\begin{displaymath}
A^\dagger(y^{(m)}) -x^{(m)} \in \hbox{ \raise-2mm\hbox{$\textstyle w-\widetilde{\lim}\atop \scriptstyle {n \to \infty}$}} \mathcal{N}(A_n) \quad (m \in \mathrm{N}).
 \end{displaymath}
Hence
\begin{displaymath}
 A^\dagger(y^{(m)})-x^{(m)} \in \mathcal{N}(A) \quad (m \in \mathrm{N}). \  (\textrm{by Lemma 3.2}).
\end{displaymath}
Then
\begin{displaymath}
x^{(m)} \in \mathcal{D}(A), \ A (A^\dagger y^{(m)}- x^{(m)})=0 \quad (m \in \mathrm{N}).
\end{displaymath}
It implies that
\begin{displaymath}
A x^{(m)}=A A^\dagger y^{(m)}\stackrel{(2.2)}{=}P_{\overline{\mathcal{R}(A)}} y^{(m)} =y^{(m)}\quad (m \in \mathrm{N})..
\end{displaymath}
\textbf{(3)}: The proof of boundedness of $ \{ x^{(m)} \}$.
\begin{eqnarray*}
\Vert x^{(m)} \Vert  & \leq \hbox{ \raise-2mm\hbox{$\textstyle \underline{\lim}\atop \scriptstyle {j \to \infty}$}}
\Vert A^\dagger_{n_j} (y^{(m)}_{n_j}) \Vert & \ \textrm{ by Lemma 2.3} \\
 & \leq  \hbox{ \raise-2mm\hbox{$\textstyle \underline{\lim}\atop \scriptstyle {j \to \infty}$}} \Vert A^\dagger_{n_j} \Vert \Vert y^{(m)}_{n_j} \Vert & \\
& \leq \sup\limits_n \Vert A^\dagger_n \Vert  \hbox{ \raise-2mm\hbox{$\textstyle \underline{\lim}\atop \scriptstyle {j \to \infty}$}} \Vert y^{(m)}_{n_j} \Vert & \ ((\textrm{by} \  (B1)) \\
 &  = \sup_n \Vert A^\dagger_n \Vert \Vert y^{(m)} \Vert & \ (\textrm{by} \ (3.12)).
\end{eqnarray*}
Taking supreme for index $ m $ on both sides yields that
\begin{displaymath}
\sup_m \Vert x^{(m)} \Vert \leq \sup_n \Vert A^\dagger_n \Vert \sup_m \Vert y^{(m)} \Vert < \infty.
\end{displaymath}
\textbf{Part II} \indent Because of Eberlein-Shmulyan theorem, the sequence $ \{ x^{(m)} \} $ contains
a weakly convergent subsequence $ \{ x^{(m_j)} \} $. Set
\begin{displaymath}
 x:= \hbox{ \raise-2mm\hbox{$\textstyle w-\lim \atop \scriptstyle {j \to \infty}$}} \  x^{(m_j)}.
\end{displaymath}
In the following, we will prove $ (x,y) \in \mathcal{G}(A) $. By Mazur theorem, for every $ j \in \mathrm{N}$, there exists a convex combination
\begin{displaymath}
\sum^{k_j}_{i=1} \alpha^{(j)}_i x^{(m_{j+i})} \ \  ( \alpha_i^{(j)} \geq 0, \sum_{i=1}^{k_j}\alpha_i^{(j)}=1 )
\end{displaymath}
such that
\begin{equation}
\Vert \sum^{k_j}_{i=1} \alpha^{(j)}_i x^{(m_{j+i})}-x \Vert \leq \frac{1}{j}.
\end{equation}
Denoting the term $ \sum\limits^{k_j}_{i=1} \alpha^{(j)}_i x^{(m_{j+i})} $ in (3.15) by  \begin{math} x_j \ (j \in \mathrm{N}) \end{math}, we rewrite (3.15) as
\begin{equation}
\Vert x_j-x \Vert \leq  \frac{1}{j}.
\end{equation}
Thus,
\begin{displaymath}
\Vert A(x_j)-y \Vert =\Vert \sum^{k_j}_{i=1} \alpha^{(j)}_i A(x^{(m_{j+i})} ) -y \Vert
\end{displaymath}.
\begin{displaymath}
=\Vert \sum^{k_j}_{i=1} \alpha^{(j)}_i y^{(m_{j+i})} -y \Vert \leq  \sum^{k_j}_{i=1} \alpha^{(j)}_i\Vert y^{(m_{j+i})} -y \Vert \end{displaymath}
\begin{displaymath}
\leq \sup_{m \geq m_j} \Vert y^{(m)}-y \Vert \leq \sup_{m \geq j} \Vert y^{(m)}-y \Vert\quad (j \in \mathrm{N})..
\end{displaymath}
 Then
 \begin{equation}
0 \leq \lim_{j \to \infty} \Vert A(x_j) -y \Vert \leq \lim_{j \to \infty} \sup_{m \geq j} \Vert y^{(m)} -y \Vert = \overline{\lim}_{j \to \infty} \Vert y^{(j)} -y \Vert =0 .
 \end{equation}
  Since \begin{math} A \end{math} is closed, we obtain from (3.16) and (3.17) that $ x \in \mathcal{D}(A)$ and $ y=Ax $. That is, \begin{math} y \in \mathcal{R}(A) \end{math}. Hence we prove that \begin{math} \mathcal{R}(A) \end{math} is closed. \qed

\section{Proof of Result (b)}
After the proof of the result (a), we obtain that, for the original operator $ A $ and its resolvent approximation setting $ \{ A_n \} $ given in Theorem 1.2, if (B1) holds, then $ A $ preserves a closed range and $ A^\dagger \in \mathcal{B}(X) $ with $ \mathcal{D}(A^\dagger) = \mathcal{R}(A)+ \mathcal{R}(A)^\perp = X $.
In the rest proof for the result (b), we only need to prove $ (C1)\Longrightarrow (B1) $ and
$ (B1) \Longrightarrow A^\dagger_n \stackrel{s}{\longrightarrow} A^\dagger $. Notice that, with Banach-Steinhaus theorem,
the former automatically holds. Thus we just need to prove the latter ($ (B1) \Longrightarrow A^\dagger_n \stackrel{s}{\longrightarrow} A^\dagger $) in the following part.
\newline \indent To prove this, we prepare a technical lemma first.
\begin{lemma}
Let $ A $, $ A_n : X \longrightarrow X , n \in \mathrm{N}$,  be bounded linear operators. Then the following two conditions are equivalent:
\newline \indent (a)  \begin{math} \mathcal{G}(A) \subseteq \hbox{ \raise-2mm\hbox{$\textstyle s-\lim \atop \scriptstyle {n \to \infty}$}} \mathcal{G}(A_n) \end{math} and \begin{math} \sup\limits_n \Vert A_n \Vert <+\infty \end{math}, where \begin{math} \Vert \cdot \Vert \end{math} is the operator norm on \begin{math} \mathcal{B}(X) \end{math};
\newline \indent (b) \begin{math} \hbox{ \raise-2mm\hbox{$\textstyle s-\lim \atop \scriptstyle {n \to \infty}$}} A_n (y) =A(y) \end{math} for every \begin{math} y \in X \end{math}.
\end{lemma}
\begin{proof}
See [10, Lemma 2.5] \qed
\end{proof}
\indent It is obvious that $ A^\dagger_n \stackrel{s}{\longrightarrow} A^\dagger$ yields from
\begin{displaymath}
\mathcal{G}(A^\dagger) \subseteq \hbox{ \raise-2mm\hbox{$\textstyle s-\lim \atop \scriptstyle {n \to \infty}$}} \mathcal{G}(A^\dagger_n)
\ \ \textrm{and} \ \  (B1): \ \sup\limits_n \Vert A^\dagger_n \Vert <+\infty
\end{displaymath}
by substituting $ A^\dagger $ and $ A^\dagger_n $ into $ A $ and $ A_n $ in Lemma 4.1 respectively.
\newline \indent Now, provided $ (B1) $ holds, under the approximation setting given in Theorem 1.2,
we are now in the position to prove $ \mathcal{G}(A^\dagger) \subseteq \hbox{ \raise-2mm\hbox{$\textstyle s-\lim \atop \scriptstyle {n \to \infty}$}} \mathcal{G}(A^\dagger_n) $.
\newline \indent Let $ (y,x) \in \mathcal{G}(A^\dagger) $, we need to construct a sequence of pairs $ \{(\xi_n, A^\dagger_n \xi_n)\} $
such that
\begin{equation}
(\xi_n, A^\dagger_n \xi_n) \stackrel{s}{\longrightarrow} (y,x).
\end{equation}
 For this construction, recalling the main idea in the proof for [12, Theorem 1.1 $(a) \Longrightarrow (b)$], we can supplement a strong convergence result for orthogonal projection sequence $ \{ P_{\mathcal{N}(A_n)} \}$ and $ \{ P_{\mathcal{R}(A_n)} \}$ in the following.
\begin{lemma}
Let $ A_n (n \in \mathrm{N}) $ and $ A $ all be defined in Theorem 1.2. If $ \{ A_n \} $ and $ A $ satisfy the resolvent consistency and (B1), then
\begin{displaymath}
\hbox{ \raise-2mm\hbox{$\textstyle s-\lim \atop \scriptstyle {n \to \infty}$}} P_{\mathcal{N}(A_n)} =P_{\mathcal{N}(A)}, \
\hbox{ \raise-2mm\hbox{$\textstyle s-\lim \atop \scriptstyle {n \to \infty}$}} P_{\mathcal{R}(A_n)}=P_{\mathcal{R}(A)}.
\end{displaymath}
\end{lemma}
\begin{proof}
\indent Recall the fact that
 \begin{displaymath}
 \hbox{ \raise-2mm\hbox{$\textstyle s-\lim \atop \scriptstyle {n \to \infty}$}} \mathcal{N}(A_n)=\mathcal{N}(A) \  \textrm{in Lemma 3.1},
 \end{displaymath}
 \begin{displaymath}
 \hbox{ \raise-2mm\hbox{$\textstyle w-\widetilde{\lim}\atop \scriptstyle {n \to \infty}$}}\mathcal{N}(A_n)
 \subset \mathcal{N}(A) \ \ \textrm{in Lemma 3.2}.
 \end{displaymath}
 Comparing the definitions of
 \begin{math}
 \hbox{ \raise-2mm\hbox{$\textstyle s-\lim \atop \scriptstyle {n \to \infty}$}} \mathcal{N}(A_n) \ \textrm{and} \ \hbox{ \raise-2mm\hbox{$\textstyle w-\widetilde{\lim}\atop \scriptstyle {n \to \infty}$}}\mathcal{N}(A_n)\end{math},
 it is obvious that
 \begin{displaymath}
 \hbox{ \raise-2mm\hbox{$\textstyle s-\lim \atop \scriptstyle {n \to \infty}$}} \mathcal{N}(A_n) \subset \hbox{ \raise-2mm\hbox{$\textstyle w-\widetilde{\lim}\atop \scriptstyle {n \to \infty}$}}\mathcal{N}(A_n).
 \end{displaymath}
Now, we have
\begin{displaymath}
\hbox{ \raise-2mm\hbox{$\textstyle s-\lim \atop \scriptstyle {n \to \infty}$}} \mathcal{N}(A_n)
\subset
\hbox{ \raise-2mm\hbox{$\textstyle w-\widetilde{\lim}\atop \scriptstyle {n \to \infty}$}}\mathcal{N}(A_n)
\subseteq
\mathcal{N}(A)=\hbox{ \raise-2mm\hbox{$\textstyle s-\lim \atop \scriptstyle {n \to \infty}$}} \mathcal{N}(A_n).
\end{displaymath}
That is,
\begin{displaymath}
\hbox{ \raise-2mm\hbox{$\textstyle s-\lim \atop \scriptstyle {n \to \infty}$}} \mathcal{N}(A_n)= \hbox{ \raise-2mm\hbox{$\textstyle w-\widetilde{\lim}\atop \scriptstyle {n \to \infty}$}} \mathcal{N}(A_n)=\mathcal{N}(A).
\end{displaymath}
With Lemma 2.2, we know
\begin{equation}
\hbox{ \raise-2mm\hbox{$\textstyle s-\lim \atop \scriptstyle {n \to \infty}$}} P_{\mathcal{N}(A_n)} =P_{\mathcal{N}(A)}.
\end{equation}
Since $ A_n $ and $ A $ are all self-adjoint, it yields that
\begin{displaymath}
\hbox{ \raise-2mm\hbox{$\textstyle s-\lim \atop \scriptstyle {n \to \infty}$}} P_{\mathcal{N}(A^*_n)} =P_{\mathcal{N}(A^*)}.
\end{displaymath}
Using identities to subtract above both sides, it follows that
\begin{displaymath}
\hbox{ \raise-2mm\hbox{$\textstyle s-\lim \atop \scriptstyle {n \to \infty}$}} P_{\mathcal{N}(A^*_n)^\perp} =P_{\mathcal{N}(A^*)^\perp}.
\end{displaymath}
Since $ \mathcal{R}(A)^\perp = \mathcal{N} (A^*)$ holds for all densely defined $ A $ on Hilbert space $ X $ (See [5, Chapter X, Proposition 1.13]), together with the fact that $ A_n (n \in \mathrm{N})$ and $ A $ possess closed ranges (Result (a)), we have
\begin{equation}
\hbox{ \raise-2mm\hbox{$\textstyle s-\lim \atop \scriptstyle {n \to \infty}$}} P_{\mathcal{R}(A_n)}=P_{\mathcal{R}(A)}.
\end{equation} \qed
\end{proof}
\noindent \textbf{Proof of Result (b)} \indent In the following, we focus on the construction of the sequence in (4.1). For any $ (x,Ax) \in \mathcal{G}(A) $, there exists a sequence $ \{ x_n \} $ such that
\begin{displaymath}
x_n \ (\in \mathcal{D}(A_n)) \stackrel{s}{\longrightarrow} x ,   \  A_n x_n \stackrel{s}{\longrightarrow} Ax \ ( n \to \infty) .  \ \ \textrm{(by (3.11))}
\end{displaymath}
Set
\begin{equation}
 z_n := P_{\mathcal{N}(A_n)^{\perp}} x_n \in \mathcal{D}(A_n) \cap \mathcal{N}(A_n)^\perp \ (\textrm{by} \ (1.3)).
 \end{equation}
 (Explanation: For $ x_n \in \mathcal{D}(A_n) = \mathcal{N}(A_n) \oplus \mathcal{C}(A_n) $, it can be uniquely represented as
\begin{displaymath}
x_n = x_{1,n} + x_{2,n}, \ \textrm{where}  \ x_{1,n} \in \mathcal{N}(A_n), \ x_{2,n} \in \mathcal{C}(A_n), \ x_{1,n} \perp x_{2,n}.
\end{displaymath}
Then
\begin{displaymath}
x_{1,n} = P_{\mathcal{N}(A_n)} x_n, x_{2,n} = x_n - x_{1,n} = x_n - P_{\mathcal{N}(A_n)} x_n = P_{\mathcal{N}(A_n)^\perp} x_n \in \mathcal{C}(A_n).)
\end{displaymath}
 Notice that,
 \begin{displaymath}
 x = A^\dagger y \in \mathcal{R}(A^\dagger) =\mathcal{D}(A) \cap \mathcal{N}(A)^\perp \subset \mathcal{N}(A)^\perp,
 \end{displaymath}
 we have
 \begin{equation}
 z_n = P_{\mathcal{N}(A_n)^{\perp}} x_n \stackrel{s}{\longrightarrow} P_{\mathcal{N}(A)^{\perp}} x = x,
 \end{equation}
 \begin{equation}
 A_n z_n \stackrel{s}{\longrightarrow} A x.
\end{equation}
Hence,
\begin{eqnarray*}
  & A_n z_n +P_{\mathcal{R}(A_n)^{\perp}}y  \stackrel{s}{\longrightarrow}  Ax+P_{\mathcal{R}(A)^{\perp}}y & \ (\textrm{by} \ (4.3) \ \textrm{and} \ (4.6)) \\
 &= P_{\mathcal{R}(A)} y + P_{\mathcal{R}(A)^{\perp}}y = y \ (\textrm{by} \ x = A^\dagger y & \ \textrm{and} \ (2.2) ),
\end{eqnarray*}
and
\begin{displaymath}
A^\dagger_n(A_n z_n +P_{\mathcal{R}(A_n)^{\perp}}y) \stackrel{(1.5)}{=} A^\dagger_n A_n z_n
 \stackrel{(2.1)}{=} P_{\overline{\mathcal{C}(A_n)}}z_n \stackrel{(4.4)}{=}z_n \stackrel{s}{\longrightarrow} x. \ (\textrm{by} \ (4.5)).
\end{displaymath}
So \begin{math} (y,x) \in \hbox{ \raise-2mm\hbox{$\textstyle s-\lim \atop \scriptstyle {n \to \infty}$}} \mathcal{G}(A^\dagger_n) \end{math}. Thus we complete the construction for (4.1). \qed
\section{Applications}
\begin{example}
Set
\begin{equation*}
 A:= - \frac{d^2}{dt^2}: \mathcal{D}(A) = H^2(0,\pi) \cap H^1_0(0,\pi) \subseteq L^2(0,\pi) \to L^2(0,\pi)
 \end{equation*}
\begin{equation*}
 A_n:= - \frac{d^2}{dt^2} - \frac{1}{n^2}: \mathcal{D}(A_n) = H^2(0,\pi) \cap H^1_0(0,\pi) \subseteq L^2(0,\pi) \to L^2(0,\pi)
 \end{equation*}
 They are all unbounded self-adjoint operators on $ L^2(0,\pi) $ (See [8, Chapter VIII.6 Example 3]).
 \newline \indent We can observe that, for  arbitrary $ y \in H^2(0,\pi) \cap H^1_0(0,\pi) $, $ A_n y \stackrel{s}{\longrightarrow} A y $ in $ L^2 $. By [15, Chapter VIII Corollary 1.6], we have
 \begin{equation*}
 R_\lambda (A_n)  \stackrel{s}{\longrightarrow} R_\lambda(A)\quad(\lambda \in \mathrm{C} \setminus \mathrm{R}), \ \textrm{and} \ (A3) \ \textrm{holds}.
 \end{equation*}
 Further, by constant coefficient variation method for ODE, we have
 \begin{equation*}
 A^\dagger y = A^{-1} y = \frac{t}{\pi} \int^\pi_0 (\pi - s ) y(s) ds - \int^t_0 (t-s) y(s) ds,
 \end{equation*}
 and
 \begin{equation*}
 A_n^\dagger y = A_n^{-1} y = \frac{n \sin \frac{t}{n}}{\sin \frac{\pi}{n}} \int^\pi_0 \sin (\frac{\pi - s }{n}) y(s) ds - n \int^t_0 \sin (\frac{t-s}{n}) y(s) ds,
 \end{equation*}
 for $ \forall y \in \mathcal{D}(A^\dagger) =\mathcal{D}(A_n^\dagger) =  \mathcal{R}(A) = \mathcal{R}(A_n)= L^2(0,\pi) $. Now it is not difficult to figure out that
 \begin{equation*}
 \Vert A_n^\dagger \Vert_{L^2 \to L^2} \leq (\frac{2\sqrt{3}}{3} +1) \pi^2, \ \forall n \  \textrm{sufficiently large}, \textrm{and} \ (B1) \ \textrm{holds},
 \end{equation*}
  \begin{equation*}
 \mathcal{D}(A^\dagger)=L^2(0,\pi), \  A^\dagger_n  \stackrel{s}\longrightarrow  A^\dagger. \ \textrm{and} \ (C1) \ \textrm{holds}.
 \end{equation*} \qed
\end{example}
\begin{example}
Set $ A = -\Delta: H^2(\mathrm{R}^3) \subseteq L^2(\mathrm{R}^3) \to L^2(\mathrm{R}^3) $,
$ A_n = -\Delta + \frac{1}{n^2}: H^2(\mathrm{R}^3) \subseteq L^2(\mathrm{R}^3) \to L^2(\mathrm{R}^3) $. They are all unbounded self-adjoint operators on $ L^2(\mathrm{R}^3) $ (see [6, Theorem 8.8]).
\newline \indent We can observe that, for arbitrary $ f \in C^\infty_0(\mathrm{R}^3) $, $ A_n f  \stackrel{s}\longrightarrow A f $. By [15, Chapter VIII Corollary 1.6], we have
 \begin{equation*}
 R_\lambda (A_n)  \stackrel{s}{\longrightarrow} R_\lambda(A) \quad (\lambda \in \mathrm{C} \setminus \mathrm{R}). \ \textrm{and} \ (A3) \ \textrm{holds}.
 \end{equation*}
 \indent Using the explicit formula for resolvent of $ -\Delta $ (See [6, Chapter 8.1]), we have
 \begin{equation*}
  A^\dagger_n = A^{-1}_n \in \mathcal{B}(L^2(\mathrm{R}^3),L^2(\mathrm{R}^3)), \ \mathcal{R}(A_n) = L^2(\mathrm{R}^3)
  \end{equation*}
  and, for arbitrary $ y \in  L^2(\mathrm{R}^3)$,
 \begin{equation*}
 (A^\dagger_n y)(s) = \int_{\mathrm{R}^3} G(s,t)y(t) dt, \ \textrm{where} \ G(s,t) = \frac{\exp(-\frac{1}{n}\vert s - t \vert )}{4\pi \vert s-t \vert}, \ s \neq t \in \mathrm{R}^3.
 \end{equation*}
 \indent Now set $ \chi := \chi_{[0,1]^3}, \Vert \chi \Vert_{L^2(\mathrm{R}^3)} =1  $. We claim that $ \lim_{n \to \infty} \Vert A^\dagger_n \chi \Vert_{L^2} =\infty $. This will yield that
 \begin{equation*}
  \lim_{n \to \infty} \Vert A^\dagger_n \chi \Vert_{L^2} \leq \overline{lim}_{n \to \infty} \Vert A^\dagger_n \Vert_{L^2 \to L^2} \leq \sup_n \Vert A^\dagger_n \Vert_{L^2 \to L^2} = + \infty.
  \end{equation*}
  And hence (B1), (C1) do not hold true.

\end{example}
Proof of the claim: Set
\begin{equation*}
u_n (s) := (A^\dagger_n \chi)(s) = \int_{\mathrm{R}^3} \frac{e^{- \frac{1}{n} \vert s-t \vert }} {4\pi \vert s-t \vert} \chi_{{[0,1]^3}} dt  = \int_{[0,1]^3} \frac{e^{- \frac{1}{n} \vert s-t \vert }} {4\pi \vert s-t \vert} dt.
 \end{equation*}
For $ \vert s \vert \geq 1 $, we have $ \vert s-t \vert \leq 2 \vert s \vert  $ and
\begin{equation*}
\exp({- \frac{1}{n} \vert s-t \vert}) \geq \exp({- \frac{2 \vert s \vert }{n}}),\  \frac{1}{{4\pi \vert s-t \vert}} \geq \frac{1}{8\pi \vert s \vert}.
\end{equation*}
Hence
\begin{equation*}
u_n (s) \geq \int_{[0,1]^3} \frac{1}{8\pi \vert s \vert} \exp({- \frac{2 \vert s \vert }{n}}) dt = \frac{1}{8\pi \vert s \vert} \exp({- \frac{2 \vert s \vert }{n}}), \ \vert s \vert \geq 1.
\end{equation*}
Then it follows that
\begin{equation*}
\Vert u_n \Vert^2_{L^2} \geq \int_{\vert s \vert \geq 1}   u^2_n (s) ds = \int_{\vert s \vert \geq 1} \frac{1}{64\pi^2 \vert s \vert^2} \exp({- \frac{4 \vert s \vert }{n}}) ds.
\end{equation*}
Using spherical polar coordinates transformation, we gain
\begin{equation*}
\Vert u_n \Vert^2_{L^2} \geq  \frac{1}{64\pi} n \exp(-\frac{4}{n})|^\infty_1  = + \infty.
 \end{equation*}
\qed
\begin{remark}
 Example 1 shows the validness of the main theorem. Example 2 tells us that the stability does not hold true simultaneously for arbitrary resolvent approximation sequence with closed range, and even a bounded scheme induced by a simple natural resolvent approximation sequence as above can be instable. Actually, we can obtain a stronger result: since $ 0 $ locates in the continuous spectrum of free Schrodinger operator (See [4, Theorem 7.17]), $ (-\Delta)^{-1} $ is unbounded, thus the (C1) does not hold and hence any bounded computational scheme $ \{ A^\dagger_n \} $ induced by resolvent approximation is instable.
\end{remark}
\begin{remark}
Above reason can also be applied to Schrodinger operators with specific potential function such that $ 0 \in \sigma_c(-\Delta + V )$, for example:
\begin{itemize}
\item Single nucleus interact with single electron (Coulomb potential):
\begin{equation*}
V(t) = - \frac{\gamma}{\vert t \vert}, \ \gamma >0, \ \mathcal{D}(-\Delta +V) =H^2(\mathrm{R}^3).
\end{equation*}
 The corresponding schrodinger operator describes the single hydrogen atom model which is propbably the most famous model in quantum mechanics. See [4, Chapter 10.2] for more information.
\item Single nucleus interact with many electrons (see [4, Chapter 11]):
\begin{eqnarray*}
H^{(N)} := - \sum^N_{j=1} \Delta_j - \sum^N_{j=1} V_{ne} (x_j) + \sum^N_{j=1} \sum^N_{j<k} V_{ee}(x_j -x_k),\\
x_j = (x_{j,1},x_{j,2},x_{j,3}) \in \mathrm{R}^3, \mathcal{D}(H^{(N)}) = H^2(\mathrm{R}^{3N}),
\end{eqnarray*}
where $ V_{ne} $ describes the interaction of one electron with the nucleus and $ V_{ee} $
describes the interaction of two electrons, and explicitly
\begin{equation*}
V_j = \frac{\gamma_j}{\vert x \vert }, \ \gamma_j > 0,\ j =ne, ee.
\end{equation*}
(Hints: $ [0,\infty) \subseteq \sigma_{ess} $ (See [4, Theorem 11.2 HVZ]), $\sigma_{p} \subseteq (\infty,0) $ (See [4, Theorem 10.4 Virial or (11.22)]) $ \Longrightarrow $ $ 0 \in \sigma_c $.)
\end{itemize}
\end{remark}
\section{Corollaries}
 In our investigation, Theorem 1.2 is the first result which extends classical Lax equivalence theorem into unbounded case with the resolvent consistency (A3). Between the resolvent consistency and the strong consistency, we know that
 \begin{proposition}
 Let $ A $ be bounded self-adjoint operator on $ X $ and $ \{ A_n \} $ a sequence of uniformly bounded self-adjoint operators, that is, $ \sup_n \Vert A_n \Vert < \infty $. Then $ (A1) \Longleftrightarrow (A3). $
 \end{proposition}
 \begin{proof}
 See [8, Chapter VIII. Problems 28].
 \end{proof}
 In this way, one can regard the resolvent consistency as a natural generalization of the strong consistency into unbounded operators. If we restrict the consideration of Theorem 1.2 in bounded case and further supplement uniform boundedness condition for approximation sequence $ \{ A_n \} $, it yields that
\begin{corollary}
 Let \begin{math} A \end{math} be bounded self-adjoint operator on Hilbert space \begin{math} X \end{math},
 \begin{math} \{ A_n \} \end{math} a sequence of uniformly bounded self-adjoint operators on \begin{math} X \end{math} with closed range \begin{math} \mathcal{R}(A_n) (n \in \mathrm{N}) \end{math}. If \begin{math} \{ A_n \} \end{math} and $ A $ satisfy the strong consistency, then
 \begin{equation*}
 \mathcal{D}(A^\dagger)=X,
 A^\dagger_n  \stackrel{s}\longrightarrow  A^\dagger  \Longleftrightarrow \sup_n \Vert A^\dagger_n \Vert <+\infty.
 \end{equation*}
 \end{corollary}
Actually a stronger version of above result which weaken the condition "uniform boundedness" into "bounded" is obtained in [12]. This version can also be deduced directly from the proof of Theorem 1.2. In particular, $ \{ A_n \} $ can be constructed in a Galerkin setting. In this way, a more specific version is obtained (in [12]) as follows:
\begin{corollary}
 Let $ A $ be a bounded self-adjoint operator on Hilbert space $ X $, $ \{ X_n \} $ a monotonically increasing and eventually dense sequence in $ X $, that is,
 \begin{equation*}
 X_n \subseteq X_{n+1} \ , \  \overline{\bigcup_{n \in \mathrm{N}} X_n} = X.
 \end{equation*}
 Set $ A_n := P_{X_n} A P_{X_n} : X \to X $, then
 \begin{equation*}
 \mathcal{D}(A^\dagger)=X,
 A^\dagger_n  \stackrel{s}\longrightarrow  A^\dagger \Longleftrightarrow \sup_n \Vert A^\dagger_n \Vert <+\infty.
 \end{equation*}
 \end{corollary}
Thus, we can regard that Theorem 1.2 provides a more unified numerical framework for approximately solving self-adjoint operator equations.
\section{Conclusion}
\indent In this work, the difficulties in essence arise from the non-closed domain of unbounded operators. To extend the Consistenncy-Stability-Convergence numerical framework into unbounded operators, the start point is to update the definition of Moore-Penrose inverse into a unbounded case. Secondly, we examine the formulation of consistenncy, stability, convergence in bounded case and provide the resolvent consistency, the core new idea in the current work, to replace the classical types of consistency. Thirdly, we follow the main idea of proof in bounded case (See [3,10,11,12,13]) to prove the main result by adjusting them with respect to non-closed domain.


%
%


\section*{Acknowledgement}
Thank Professor Nailin Du for indicating the main theorectical result and the guidance throughout the whole research process, and Doctor Jochen Glueck for the help in spectral theory.
\section*{Conflict of interest}
 The author declares that he has no conflict of interest with someone who concerns this work.

\bibliographystyle{spbasic}      

\end{document}